\documentclass[11pt,a4paper]{article}
\usepackage[english]{babel}
\usepackage{hyperref}
\usepackage[utf8]{inputenc}

\usepackage{color}

\usepackage{fontenc}

\usepackage{amssymb, amsmath, amsbsy}
\usepackage{array}
\usepackage{amsthm}

\usepackage{fancyhdr}

\usepackage{nccmath}

\usepackage{geometry}

\usepackage{graphicx} 
\usepackage{subfigure} 
\usepackage{float} 

\newcommand{\<}{\left<}
\renewcommand{\>}{\right>}

\newcommand{\com}[1]{\opt{draft}{\textcolor{red}{
$\LHD$ #1 $\RHD$\marginpar{\textcolor{red}{$\begin{lema}acksquare$}}}}}

\newcommand{\comb}[1]{\opt{draft}{\textcolor{blue}{
$\LHD$ #1 $\RHD$\marginpar{\textcolor{blue}{$\begin{lema}acksquare$}}}}}

\newenvironment{demo}{{\bf Proof }}
{\qed \\}

\newcommand{\arccoth}{{\rm arccoth}}

\newcommand{\re}{\mathbb R}

\renewcommand{\(}{\left(}
\newcommand{\lb}{\label}
\newcommand{\nn}{\nonumber}

\newcommand{\fracc}{\displaystyle\frac}
\newcommand{\ds}{\displaystyle}
\renewcommand{\)}{\right)}
\newcommand{\flecha}{\longrightarrow}

\newcommand{\eps}{\ensuremath{\varepsilon}}
\def\a{\alpha}
\def\p{\varphi}

\def\vec{\overrightarrow}

\def\k{\frak k}

\def\s{{\rm s_\lambda}}
\def\c{{\rm c_\lambda}}

\def\co{{\rm co_\lambda}}

\def\ta{{\rm ta_\lambda}}
\def\co{{\rm co_\lambda}}

\newcommand{\bde}{\begin{defi}}
\newcommand{\ede}{\end{defi}}

\newcommand{\be}{\begin{enumerate}}
\newcommand{\ee}{\end{enumerate}}

\newcommand{\ba}{\begin{array}}
\newcommand{\ea}{\end{array}}

\def\oMu{{\overline{M}_\lambda^{n+1}}}
\def\oMd{{\overline{M}_\lambda^{2}}}

\def\oA{{\overline A}}
\def\orO{{\overline O}}

\def\hM{{\widehat M }}

\def\hB{{\widehat {\rm B}}}

\def\hA{{\widehat {A}}}
\def\hC{{\widehat {C}}}

\def\raiz{\sqrt{\tau}}

\newtheorem{defi}{\hspace{12pt} Definition}
\newtheorem{teor}{\hspace{12pt} Theorem}
\newtheorem{prop}[teor]{\hspace{12pt} Proposition}

\newtheorem{lema}[teor]{\hspace{12pt} Lemma}
\newtheorem*{lema*}{\hspace{12pt} Lemma}
\newtheorem*{teor*}{\hspace{12pt} Theorem}

\def\vec{\overrightarrow}

\def\k{\frak k}

\def\s{{\rm s_\lambda}}
\def\c{{\rm c_\lambda}}

\def\co{{\rm co_\lambda}}

\def\ta{{\rm ta_\lambda}}
\def\co{{\rm co_\lambda}}

\def\h{{\kappa}}

\newcommand{\ben}{\begin{enumerate}}
\newcommand{\een}{\end{enumerate}}
\newcommand{\bi}{\begin{itemize}}
\newcommand{\ei}{\end{itemize}}
\newcommand{\bec}{\begin{equation}}
\newcommand{\eec}{\end{equation}}
\newcommand{\beca}{\begin{equation*}}
\newcommand{\eeca}{\end{equation*}}
\newcommand{\bal}{\begin{align}}
\newcommand{\aal}{\end{align}}
\newcommand{\bala}{\begin{align*}}
\newcommand{\aala}{\end{align*}}

\begin{document}

\title{Inscribed and circumscribed radius of $\kappa$-convex hypersurfaces in Hadamard manifolds \footnotetext{The research is partially supported by grant PID2019-105019GB-C21 funded by MCIN/AEI/ 10.13039/501100011033 and by \lq\lq ERDF A way of makimg Europe\rq\rq , and by the grant AICO 2021 21/378.01/1 funded by the Generalitat Valenciana.}
}
\author{Alexander Borisenko and Vicente Miquel}

\date{}

\maketitle

\begin{abstract}

Let $P$ be a convex polygon in a Hadamard surface $M$ with curvature $K$ satisfying $-k_2^2 \ge K \ge -k_1^2$.  We give an upper bound of the circumradius of $P$ in terms of a lower bound of the curvature of $P$ at its vertices.
\end{abstract}


\section{Introduction and main result}

Let $M$ be a complete  $m$-dimensional Riemannian  manifold. Given a domain $\Omega \subset M$, an inscribed ball or inball is a ball in $M$ contained in $\Omega$ with maximum radius, which is called the inradius of $\Omega$, and we shall denote it by $r$.  A circumscribed ball of $\Omega$ is a ball in $M$ containing $\Omega$ with minimum radius, which is called the circumradius and we shall denote it by $R$.

In the book of Blaschke \cite{Bl49} it is proved that
{\it  if $\Gamma$ is a closed convex regular curve in the Euclidean plane that bounds a compact convex region $\Omega$ and the curvature $\h$ of $\Gamma$ is bounded from below by some constant $\h_0>0$, then, the circumradius $R$ of $\Omega$ is bounded from above by  $\fracc1{\h_0}$.}

This result was extended by H. Karcher (\cite{Ka68}) for the other space forms. Further developments of related Blaschke theorems has been done by obtaining conditions under which a convex set in $\re^n$ can be included  in other (\cite{Ra74,De79,BS89}).

In all these theorems, the hypothesis of strong convexity ($\h\ge\h_0>0$) is necessary, the theorem is not true for $\h\ge 0$. Then  it cannot be applied to closed convex polygons. In \cite{BM23} we used two definitions of curvature at the vertex of a polygon which allowed us to obtain a version of the Blaschke Theorem for polygons in space forms. Here we use the same definitions to obtain a less precise upper bound of the circumradius $R$ for polygons in Hadamard surfaces with curvature $-k_2^2 \ge K \ge -k_1^2$.

In order to do that, first we use comparison theorems to bound the inradius $r$ of a domain. Then, we apply Th. 3.1 in \cite{BM02} to get an upper bound of $R$. Then we use this to obtain the upper bound of $R$ for polygons. This is done by a  combination of the comparison theorems mentioned above with the comparison of the angles given by Toponogov Theorem, which has to be done carefully because these inequalities  go in opposite sense.

Let us recall some known definitions and properties  before giving precise statements.

We shall work on an $n$-dimensional Hadamard manifold $M^n$, that is, a simply connected  complete $n$-dimensional Riemannian manifold with sectional curvature bounded from above by a constant $-k_2^2\le0$). In such a manifold, we define

\begin{defi} Let $\lambda >0$, an orientable smooth ($C^2$ or more) hypersurface $L$ of a Hadamard manifold $M^n$ is called $\lambda$-convex if there is a suitable selection of the unit normal vector of $L$ such that the normal curvatures $k_N$ of $L$ satisfy $k_N\ge\lambda$. 
\end{defi}

\begin{defi}
A domain $\Omega \subset M$ is called convex if each shortest path with endpoints in $\Omega$ lies in $\Omega$.
\end{defi}
If $\Omega$ is convex in $M$, then $\partial \Omega$ is a topological embedded hypersurface which is smooth except for a set of zero measure.

\begin{defi} A convex domain $\Omega$ of $M$ is $\lambda$-convex if for every point $p\in\partial \Omega$ there is a smooth $\lambda$-convex hypersurface $L$ through $p$ such that there is a neighborhood of $p$ in $\Omega$ contained in he convex side of $L$ (that is, the side in $M$ where the unit normal vector to $L$ points). 
\end{defi}
It is known that if $\partial \Omega $ is smooth, then $\Omega $ is $\lambda$-convex if and only if $\partial \Omega$ is $\lambda$-convex.

\begin{defi}
A topological immersion $f:P\flecha M$ is called locally convex at a point $x\in P$ if $x$ has a neighborhood $U$ in $P$ such that $f(U)$ is part of the boundary of a convex domain of $M$.
\end{defi}

In  \cite{Al77,Bo00}), it was proved that 

\begin{teor*}[\cite{Al77,Bo00}]
If $P$ is a compact orientable locally convex and immersed hypersurface of dimension $n\ge 2$ in a Hadamard manifold $M$, then $P$ is embedded, homeomorphic to the sphere and it is the boundary of a convex set $\Omega$.
\end{teor*} 

For $n=1$, the above theorem is not true even when the immersion is $C^2$. In this paper, even for dimension $1$ we shall consider only compact embedded convex curves which are the boundary of a convex domain $\Omega\subset M$, that is, for any dimension, we shall adopt the following

\begin{defi}
A compact orientable locally convex hypersurface $P$ in $M$ is $\lambda$-convex if it is the boundary of a $\lambda$-convex domain $\Omega$.

By the inradius and the circumradius of $P$ we understand the inradius and circumradius of $\Omega$.
\end{defi} 

We remark that, with this general definition, $P$ can be $\lambda$-convex and, at the same time, contain conical or rudge points, or points where $P$ is only $C^1$, where it allowed to say that normal curvature is infinite in some directions.

In the next sections we shall prove the following:

\begin{teor}\lb{te2}
Let $M$ be a Hadamard manifold with sectional curvature $K$ satisfying $-k_2^2\ge K \ge - k_1^2$, $k_2, k_1>0$. Let $P$ be a compact hypersurface of $M$ such that
\bec\lb{hipt1}
P \text{ is } k_2 \coth(k_2\rho)-\text{convex } \text{ and } \ k_2 \coth(k_2\rho) \ge k_1,
\eec
 then 
 \begin{align}\lb{1-1}
 & r\le \frac{1}{k_1} \arccoth\(\frac{k_2}{k_1} \coth(k_2\rho)\) \ \text{ and }\ 
  R\le \frac{1}{k_1} \arccoth\(\frac{k_2}{k_1} \coth(k_2\rho)\) + k_1 \ln 2. 
 \end{align} 
\end{teor}

\begin{defi}
In a surface $M$, let $P$ be a polygon. If $A$ is a vertex of the polygon, $\a$ the interior angle at $A$ and $\ell_1$, $\ell_2$ the lengths of the sides of $P$ that meet at vertex $A$, then the curvature of $P$ at $A$ is defined by:
\bec \lb{kA}
\h_A = \frac{2\ (\pi-\a)}{\ell_1+\ell_2}. 
\eec
\end{defi}

\begin{teor}\lb{te3}
Let $M$ be a Hadamard surface with Gauss curvature $K$ satisfying $-k_2^2> K \ge - k_1^2$, $k_2,k_1>0$. Let $P$ be a polygon with 
 sides of lengths $\ell_i$ and  vertices $A_i$. 
If $\h_{A_i} \ge \frac{\pi}{2} k_1 \coth(k_1\rho)$ and $\coth(k_2\rho) \ge \frac{k_1}{k_2}$, then the inradius $r$ and the circumradius $R$ of $P$ satisfy 
\begin{align}
& k_1 \coth(k_1 r) \ge  k_2 \coth(k_2\rho), \quad \text{ that is:} \lb{rupb} \\
&r\le \frac1{k_1} \arccoth\(\frac{k_2}{k_1} \coth(k_2 \rho)\) \ \text{ and }\  R \le \frac1{k_1} \arccoth\(\frac{k_2}{k_1} \coth(k_2 \rho)\) + k_1 \ln 2.
\end{align}
\end{teor}

\section{Proof of Theorem \ref{te2}}

We shall use the following result:

\begin{lema}[\cite{BM02},Th. 3.1]\lb{had} Let $M$ be a Hadamard manifold
with sectional curvature $K$ satisfying $0\ge K \ge
-k_1^2$. 
 If
$\Omega$ is a  compact 
$k_1$-convex domain, then 
\begin{align}\lb{bm02}
R - r
\leq k_1 \ln\frac{(1+\raiz)^2}{1+\tau} < k_1 \ln 2,
\end{align}
where $\tau = \tanh(k_1r/2)$.
Moreover this bound is sharp.
\end{lema}


Let $S$ be the geodesic sphere of $M$ which is the boundary of an inball of $P$. From standard comparison theory (see \cite{Pe06}) we have that the normal curvature of $S$ at any point satisfies 
\bec\lb{kns}
k_2 \coth(k_2 r)\le k_N^S\le k_1 \coth(k_1 r).
\eec
Let $Q_0 \in S\cap P$.
 Since $P=\partial\Omega$, which is $k_2 \coth(k_2 \rho)$-convex, there is a smooth $k_2 \coth(k_2 \rho)$-convex hypersurface $L$ through $Q$ leaving a neighbourhood of $Q$ in $\Omega$ (and then in $S$) in the convex side of $L$, then $S$ and $L$ are tangent at $Q$ and  \bec\lb{knl}
 k_N^S(Q_0) \ge k_N^L(Q_0) \ge k_2 \coth(k_2 \rho),
 \eec
  where $k_N^L(Q_0)$ is the normal curvature of $L$ at $Q_0$. From \eqref{kns} and \eqref{knl}  we obtain that 
  \bec\lb{rbound}
   k_2 \coth(k_2 \rho) \le k_N^L(Q_0)\le k_1 \coth(k_1 r),
   \eec
from which we obtain the first inequality of \eqref{1-1}. Now, by using the hypothesis \eqref{hipt1} we have  that $P$ is   $k_1$-convex, then we can apply Lemma  \ref{had} to obtain the second inequality of \eqref{1-1}. \qed

\section{Proof of Theorem \ref{te3}}

Let $\ell_i$, $\ell_{i+1}$ be the lengths of the sides having $A_i$ as a common vertex. As in the case of constant curvature (see \cite{BM23}), we consider the segments of circles $C_i$ from $A_{i-1}$ to $A_i$ of radius $\rho_i$ and center $O_i$ and $C_{i+1}$ from $A_i$ to $A_{i+1}$ of radius $\rho_{i+1}$ and center $O_{i+1}$. Now, on the Hyperbolic space $\oMd$ of constant sectional curvature $\lambda = -k_1^2$, we consider the geodesic triangles $\orO_i \oA_{i-1} \oA_i$, $\orO_{i+1} \oA_{i} \oA_{i+1}$ with sides with the same lengths than the the corresponding triangles $O_i A_{i-1} A_i$, $O_{i+1} A_{i} A_{i+1}$ in $M$. From the Toponogov comparison theorem on the angles of a triangle we have that interior angles of the triangles in $M$ are bigger than the corresponding ones in $\oMd$.

We want that the curve obtained by the union os the segments of circle $C_i$ be convex, and this will happen if and only if the angle $\widehat{O_{i+1}A_iO_i} \in [0, \pi]$, but this occurs if and only if 
\bec\lb{cco1}
\widehat{A_{i-1}A_iO_i} + \widehat{O_{i+1}A_iA_{i+1}} \ge \widehat{A_{i-1}A_iA_{i+1}}.
\eec
 If $\delta_i = \frac{\pi}2 - \widehat{A_{i-1}A_iO_i}$ and $\delta_{i+1} = \frac{\pi}2 - \widehat{O_{i+1}A_iA_{i+1}}$, we can use the definition \eqref{kA} to write the inequality \eqref{cco1} under the form:
 \bec\lb{cco2}
 \delta_i + \delta_{i+1} \le \pi - \widehat{A_{i-1}A_iA_{i+1}} = \kappa_{A_i}(\ell_i+\ell_{i+1})/2.
 \eec
  On the other hand we have $\delta_i =  \pi/2- \widehat{A_{i-1}A_iO_i} \le \pi/2 -  \widehat{\oA_{i-1}\oA_i\orO_i} =: \overline \delta_i$ , $\delta_{i+1}=  \pi/2-\widehat{O_{i+1}A_iA_{i+1}} ) \le \pi/2-  \widehat{\orO_{i+1}\oA_i\oA_{i+1}}=:\overline\delta_{i+1}$. Then, the inequality \eqref{cco2} is satisfied if 
 \bec\lb{cco3}
 \overline\delta_i + \overline\delta_{i+1} \le  \kappa_{A_i}(\ell_i+\ell_{i+1})/2,
 \eec 
 is satisfied. Now we are going to check that the hypothesis on the lower bound of $\h_{A_i}$ implies \eqref{cco3}.

From hyperbolic trigonometry applied to the triangles $\orO_i \oA_{i-1} \oA_i$, $\orO_{i+1} \oA_{i} \oA_{i+1}$ one has that $\tanh(k_1\ell_i/2) = \tanh(k_1\rho_i)\sin(\overline \delta_i)$ for $i$ and for $i+1$.

Since $\overline \delta_i \in[0,\pi/2]$, $\overline \delta_i \le \pi/2 \sin(\overline \delta_i)$. Moreover, we shall take $\rho_i = \rho$. Then

\begin{align}
&\overline \delta_i + \overline \delta_{i+1} \le \frac{\pi}{2} (\sin \overline\delta_i+ \sin \overline\delta_{i+1}) = \frac{\pi/2}{ \tanh(k_1\rho) }  ( \tanh(k_1\ell_i/2)+  \tanh(k_1\ell_{i+1}/2))\nn \\
&  = \frac{\pi/2}{\tanh(k_1\rho) }  \frac{ \tanh(k_1\ell_i/2)+  \tanh(k_1\ell_{i+1}/2)}{(\ell_i + \ell_{i+1})/2} (\ell_i + \ell_{i+1})/2 \nn \\
&= \frac{ \tanh(k_1\ell_i/2)+ \tanh(k_1\ell_{i+1}/2)}{(\ell_i + \ell_{i+1})/2} \frac{\pi/2}{ \tanh(k_1\rho) \kappa_{A_i}} \kappa_{A_i} (\ell_i + \ell_{i+1})/2 \nn \\
& \le \frac{ \tanh(k_1\ell_i/2)+ \tanh(k_1\ell_{i+1}/2)}{k_1 (\ell_i + \ell_{i+1})/2}  \kappa_{A_i}  (\ell_i + \ell_{i+1})/2 \le \kappa_{A_i}  (\ell_i + \ell_{i+1})/2\lb{bpAi}
\end{align}
which is the desired condition.

Now, let us take a parallel $C_\eps$ to $C$ at distance $\eps$. $C_\eps$ is the union of segments of circles $C'_i$ of radius $\rho+\eps$ and center $O_i$ and circles of radius $\eps$ centered at $A_i$, then it is $C^{1,1}$ and its normal curvature  $k_N$ satisfies $k_2 \coth(k_2 \eps) \le k_N \le k_1 \coth(k_1 \eps)$ at points of $C_\eps$ at distance $\eps$ from the vertices $A_i$ and $k_2 \coth(k_2 \rho+\eps) \le k_N \le k_1 \coth(k_1 \rho+\eps)$ for others. Then, for every $\eps$, $k_2 \coth(k_2 \rho_i+\eps) \le k_N$, and we can apply Theorem \ref{te2} to conclude  $k_1 \coth(k_1 r) \ge k_2 \coth(k_2 (\rho+\eps))$ for $C_\eps$, then for $P$ because the domain bounded by $P$ is included in the domain bounded by $C_\eps$. Taking $\eps\to 0$, we obtain $k_1 \coth(k_1 r) \ge k_2 \coth(k_2 \rho)$, which is \eqref{rupb}.  \qed

\section{Remarks}

\subsection{Bounds for Theorem\ref{te2} in terms of  only $k_1$}

We have stated Theorem \ref{te2} under a form that  has a direct application for the proof of Theorem \ref{te3}. This form implies some restrictions for the number $\rho$ ($\rho \le \frac1{k_2} \arccoth(\frac{k_1}{k_2})$), which appears in the lower bound of the normal curvatures of the hypersurface $P$. The same arguments than for the proof of Theorem \ref{te2} allow to obtain  another upper bound for $r$ with hypotheses that  impose no restriction on $\rho$. The statement of this other result is:

\begin{teor*}{\bf 1'}
Let $M$ be a Hadamard manifold with sectional curvature $K$ satisfying $0\ge K \ge - k_1^2$, $k_1>0$. Let $P$ be a compact $k_1 \coth(k_1\rho)$-convex  hypersurface of $M$, then 
 \begin{align}
 & r\le \rho \text{ and }
  R\le \rho + k_1 \ln 2. 
 \end{align} 
\end{teor*}

\subsection{The theorems when $k_2=0$}

In this case, in the hypotheses, the lower bounds for $k_N$ or $k_{a_i}$ will be  $1/\rho$ instead of $k_2 \coth(k_2\rho)$, and $\rho$ has to satisfy $1/\rho \ge k_1$. Then, the statement of the theorems, under similar proof that the ones given before, is

\begin{teor*}{\bf 1''}
Let $M$ be a Hadamard manifold with sectional curvature $K$ satisfying $0\ge K \ge - k_1^2$, $k_1>0$. Let $P$ be a compact hypersurface of $M$ such that
\bec
P \text{ is } \frac1\rho-\text{convex } \text{ and } \frac1\rho \ge k_1,
\eec
 then 
 \begin{align}
 & r\le \frac{1}{k_1} \arccoth\(\frac{1}{k_1\ \rho}\) \quad \text{ and } \quad 
  R\le \frac{1}{k_1} \arccoth\(\frac{1}{k_1\ \rho}\) + k_1 \ln 2. 
 \end{align} 
\end{teor*}

\begin{teor*}{\bf 2'}
Let $M$ be a Hadamard surface with Gauss curvature $K$ satisfying $0 > K \ge - k_1^2$, $k_1>0$. Let $P$ be a polygon with $n$ sides of lengths $\ell_i$ and  vertices $A_i$. 
If $\h_{A_i} \ge \frac{\pi}{2} k_1 \coth(k_1\rho)$ and $\frac1\rho \ge k_1$, then the inradius $r$ and the circumradius $R$ of $P$ satisfy 
\begin{align*}
& k_1 \coth(k_1 r) \ge  \frac{1}{\rho}, \text{ that is:} \quad
r\le \frac1{k_1} \arccoth\(\frac{1}{k_1 \rho} \) \ \text{ and }\  R \le \frac1{k_1}\arccoth\(\frac{1}{k_1 \rho}\) + k_1 \ln 2.
\end{align*}
\end{teor*}

\subsection{Theorem \ref{te3} with other definition of $k_A$.}

In \cite{BM23} another definition of the curvature of a polygon in a surface of constant sectional curvature is done. It coincides with \eqref{kA}, but in spaces of constant sectional curvature $-k_1^2$ it takes the form

\bec \lb{kAnc}
\h_A = \frac{(\pi-\a)}{\frac1{k_1}\tanh( k_1 \ell_1/2)+\frac1{k_1}\tanh( k_1 \ell_2/2)}. 
\eec

Taking this definition for surfaces with $0\ge K \ge - k_1^2$ clould
seem as natural as taking definition \eqref{kA}. With the computations that we have done, the result will be again Theorem \ref{te3}. The reason is that in the last inequality of \eqref{bpAi}, we have bounded $\ds\frac{ \tanh(k_1\ell_i/2)+ \tanh(k_1\ell_{i+1}/2)}{k_1 (\ell_i + \ell_{i+1})/2}$ by $1/k_1$, but, with  the new definition of $\kappa_{A_i}$, inequality \eqref{cco3} changes to $\overline\delta_i + \overline\delta_{i+1} \le  \kappa_{A_i}(\frac1{k_1}(\tanh(\ell_i/2)+\tanh(\ell_{i+1})/2)$, and the above quotient is $1/k_1$, the same value of the bound we have taken.

\subsection{The theorems when $k_2=k_1$}

For Theorem \ref{te2}, the hypothesis will be only $P$ is $k_21\coth(k_1 \rho)$-convex and the thesis will be  $r\le \rho$ and $R\le \rho+k_1 \ln 2$. If we compare this with the corresponding Theorem by Karcher when $P$ is $C^2$, where, with the same hypothesis, we obtain $R\le \rho$ we see that our bound is not the best one: we bounded $r$ for which should be a bound of $R$. Then or bound is far from being the best one.

Similar remarks on Theorem \ref{te3}: The hypothesis will be \lq\lq the same \rq\rq, and the conclusion will be $r\le \rho$. In \cite{BM23} we obtained, with the same hypothesis, $R\le \rho$, again the same difference than with Theorem \ref{te2}.

\bibliographystyle{alpha}

\begin{thebibliography}{99}

\bibitem {Al77} Stephanie Alexander, Locally convex hypersurfaces of negatively curved spaces
{\it
Proc. Amer. Math. Soc.} 64 (1977), no. 2, 321–325.

\bibitem {Bl49}  Wilhelm  Blaschke; {\it Kreis und Kugel}. Chelsea Publishing Co., New York, 1949. x+169 pp. (Photo-offset reprint of the edition of 1916 [Veit, Leipzig]) 

\bibitem {Bo00} Borisenko, A. A., On locally convex hypersurfaces in Hadamard manifolds.
{\it Math. Notes} 67 (2000), no. 3-4, 425–432.

\bibitem{BM23} A. A. Borisenko and V. Miquel "A discrete Blaschke Theorem for convex polygons in 2-dimensional space forms" arXiv:2305.07566  

\bibitem{BM02} A. A. Borisenko and V. Miquel "Comparison theorems on convex hypersurfaces in Hadamard manifolds" {\it Ann. Global Anal. Geom.} 21 (2002), no. 2, 191–202.

\bibitem{BS89} Jeff Brooks and John B. Strantzen; Blaschke’s rolling theorem in $\re^n$.
{\it Mem. Amer. Math. Soc.} 80 (1989), no. 405, vi+101 pp.

\bibitem{De79} José A. Delgado; Blaschke's theorem for convex hypersurfaces. {\it J. Differential Geometry} 14 (1979), no. 4, 489–496 (1981). 


\bibitem{Ka68}  Hermann Karcher; Umkreise und Inkreise konvexer Kurven in der sphärischen und der hyperbolischen Geometrie. {\it Math. Ann.} 177 (1968), 122–132.


\bibitem{Pe06} Peter Petersen, {\it Riemannian geometry},
Grad. Texts in Math., 171
Springer, New York, 2006.

\bibitem{Ra74} Jeffrey Rauch; An inclusion theorem for ovaloids with comparable second fundamental forms. {\it J. Differential Geometry} 9 (1974), 501–505. 


\end{thebibliography}

{ Alexander Borisenko, \\
{B. Verkin Institute for Low Temperature, Physics and Engineering of the National Academy of Sciences of Ukraine\\
Kharkiv, Ukraine
}\\ and \\
 Department of Mathematics \\
University of Valencia \\
46100-Burjassot (Valencia), Spain}

{aborisenk@gmail.com}\\

{ Vicente Miquel \\
Department of Mathematics \\
University of Valencia \\
46100-Burjassot (Valencia), Spain
}

{vicente.f.miquel@uv.es}

\end{document}